\documentclass[11pt]{amsart}

\usepackage{amsfonts,amssymb,amscd,amsmath,latexsym,amsbsy,amsthm}
\usepackage{url}
\usepackage{amssymb}
\usepackage{amsfonts}
\usepackage{latexsym}

\theoremstyle{plain}
\newtheorem{theorem}{Theorem}

\newtheorem{lemma}[theorem]{Lemma}

\newtheorem{proposition}[theorem]{Proposition}

\newtheorem{corollary}[theorem]{Corollary}

\newtheorem{conjecture}[theorem]{Conjecture}

\theoremstyle{definition}

\theoremstyle{remark}

\numberwithin{theorem}{section}
\numberwithin{equation}{section}

\newcommand{\Tr}{\text{Tr}}

\newcommand{\g}{\mathfrak{g}}
\newcommand{\h}{\mathfrak{h}}
\newcommand{\n}{\mathfrak{n}}

\newcommand{\ben}{\begin{enumerate}}
\newcommand{\een}{\end{enumerate}}
\newcommand{\ad}{{\text{Ad}}}

\newcommand{\solu}[1]{\begin{sol}{\bf (\ref{#1})}}

\begin{document}

\title{On central extensions of preprojective algebras}

\author{Pavel Etingof}
\email{etingof@math.mit.edu}

\author{Fr\'ed\'eric Latour}
\email{flatour@mailbox.syr.edu}

\author{Eric Rains}
\email{rains@math.ucdavis.edu}
\maketitle

\centerline{Dedicated to the 70th birthday of Ernest Vinberg}

\section{Introduction}

Let $Q$ be a quiver of ADE type. Let $\overline{Q}$ be the double
of $Q$, and $P$ the path algebra of $\overline{Q}$ over $\Bbb
C$. The paper \cite{ER} attaches to $Q$ a centrally extended preprojective algebra
$A=A^\mu$, which is the quotient of $P[z]$
by the relation $\sum_{a\in Q}[a,a^*]=z(\sum \mu_ie_i)$,
where $\mu=(\mu_i)$ is a regular weight (for the root system attached
to $Q$), and $e_i$ are the vertex
idempotents in $P$.\footnote{This algebra is denoted in \cite{ER}
by $\Pi_0^\mu$.} It is shown in
\cite{ER} that the algebra $A$ has nicer properties
than the ordinary Gelfand-Ponomarev
preprojective algebra $A_0=A/(z)$ of $Q$; in particular,
the deformed version $A(\lambda)$ of $A=A(0)$
is flat, while this is not the case for $A_0$.
The paper \cite{ER} also shows that $A$ is a Frobenius algebra,
and computes the Hilbert series of $A$. Finally, \cite{ER}
links the algebra $A$ with cyclotomic Hecke algebras of complex
reflection groups of rank 2.

The goal of this paper is to continue to study the rich structure of
the algebra $A$. In particular, we show that
for generic $\mu$ (and specifically for $\mu=\rho$) the algebra $A$ has a unique
trace functional, and compute the structure of the center $Z$
of $A$ and the trace space $A/[A,A]$. Namely, it turns out
that $Z$ and $A/[A,A]$ are dual to each other under the trace
form, and the dimension of the homogeneous subspace
$(A/[A,A])[2p]$ equals the number of positive roots for $Q$ of
height $p+1$.

We also show that the elements $z^s(\sum \phi_i e_i)$
span $A/[A,A]$, and determine when such an element
maps to zero in $A/[A,A]$ (i.e. sits in $[A,A]$).
The answer is given in terms
of the structure of the maximal nilpotent subalgebra $\n$
of the simple Lie algebra $\g$ attached to $Q$, which demystifies
the equality between $\dim(A/[A,A])[2p]$ and the number of
positive roots for $Q$ of height $p+1$.

{\bf Acknowledgments.} P.E. thanks George Lusztig for a useful
discussion. The work of P.E. was partially supported by the NSF
grant DMS-0504847 and the CRDF grant RM1-2545-MO-03. F.L. thanks
the IH\'ES in Bures-sur-Yvette, France. E.R. was supported in part
by NSF grant DMS-0401387.

\section{Preliminaries and some results}

\subsection{Preliminaries}

We recall some definitions and notation from \cite{ER}.

Let $Q$ be a quiver of ADE type.
Let $I$ be the set of vertices of $Q$, and $r=|I|$.

Consider the root system ${\mathcal R}$ attached to $Q$.
Let $\alpha_j, j\in I$, be the simple roots.
Let $\omega_j,j\in I$, be the fundamental weights.
Let $\rho=\sum \omega_i$. If $\alpha$ is a positive root, then
the height of $\alpha$ is the number of simple roots occuring
in the decomposition of $\alpha$; it equals to the inner product
$(\rho,\alpha)$. Let $h$ be the Coxeter number of ${\mathcal R}$.
Let $N$ be the number of positive roots in ${\mathcal R}$.
Recall that $N=hr/2$.

Let $\g$ be the simple Lie algebra whose Dynkin diagram is $Q$.
Fix a polarization $\g=\n_+\oplus \h\oplus \n_-$, where $\n_\pm$
are the nilpotent subalgebras, and $\h$ the Cartan subalgebra.
For brevity we will denote $\n_-$ by $\n$.
The Lie algebra $\n$ is generated by elements $F_i, i\in I$,
subject to the Serre relations.

Let $R$ be the algebra of complex-valued functions on $I$, and
$e_i,i\in I$, be the primitive idempotents of this algebra.
Let $\overline{Q}$ be the double of $Q$.
Let $V$ be the $R$-bimodule spanned by the edges
of $\overline{Q}$.
Let $P=T_RV$ be the path algebra of the doubled quiver
$\overline{Q}$ (the tensor algebra over $R$ of the bimodule
$V$). Let $\mu=\sum_{i\in I}\mu_i\omega_i\in \h^*$ be a regular
weight (i.e. the inner product $(\mu,\alpha)\ne 0$ for any root
$\alpha\in {\mathcal R}$).
Define the centrally extended preprojective algebra $A=A^\mu$
of $Q$, which is the quotient of $P[z]$
(where $z$ is a central variable) by the
relation
$$
\sum_{a\in Q}[a,a^*]=z(\sum_{i\in I}\mu_ie_i).
$$
Note that if $\mu=\rho$ then this relation takes an especially
simple form
$$
\sum_{a\in Q}[a,a^*]=z.
$$
Also, let $A_0:=A/(z)$ be the usual preprojective algebra of $Q$
(it is the quotient of $P$ by the relation $\sum_{a\in Q}[a,a^*]=0$.)

Define the deformed centrally extended preprojective algebra
$A(\lambda)=A^\mu(\lambda)$ to be the quotient of the
path algebra $P[z]$ by the defining relation
$$
\sum_{a\in Q} [a,a^*]=\sum_{i\in I}(\mu_iz+\lambda_i)e_i,
$$
where $\lambda=\sum_{i\in I} \lambda_i\omega_i\in \h^*$ is a weight.
This algebra carries a natural filtration, given by
$\deg(R)=0$, $\deg(a)=\deg(a^*)=1$, $\deg(z)=2$.
It is shown in \cite{ER} that $A(\lambda)$ is a flat
deformation of $A(0)=A$, i.e., ${\rm gr}(A(\lambda))=A(0)$.

It is clear that the algebras $A_0$ and $A(\lambda)$ are
independent on the orientation of $Q$, up to an isomorphism.

\subsection{The trace function on $A$}

 From now on we assume that $\mu$ is a fixed generic weight,
or $\mu=\rho$.

Recall from \cite{ER} that $A$ is a finite dimensional
$\Bbb Z_+$-graded Frobenius algebra, with socle in
degree $2(h-2)$, with basis $z^{h-2}e_i$.\footnote{Note that the elements
$z^{h-2}e_i$ may vanish for special regular $\mu$.}

\begin{proposition}\label{trace}
(i) There exists a unique up to scaling trace $\Tr: A\to \Bbb C$ of
degree $2(h-2)$, i.e. a nonzero linear functional such that
$\Tr(xy)=\Tr(yx)$.

(ii) The form $(x,y):=\Tr(xy)$ is nondegenerate.
\end{proposition}

\begin{proof} Clearly, we may assume that $Q$ has at least two
vertices. Recall that the degree 1 component $A[1]$ of $A$ is spanned by
edges $a$ of the doubled quiver $\overline Q$.
Also, $A[2(h-2)-1]$ is
spanned by elements of the form $z^{h-3}b$, where $b$
is an edge of $\overline{Q}$. Indeed, it follows from \cite{ER},
Section 4, that this is true for $\mu=\rho$, hence it is true for
generic $\mu$ by deformation argument.

Since $A$ is a Frobenius algebra, for every edge $a$ we have
$z^{h-3}aa^*=c_a z^{h-2}e_{head(a)}$,
where $c_a$ is a nonzero
number.

But $[A,A][2(h-2)]=[A[1],A[2(h-2)-1]]$, so it is the span of
$z^{h-3}[a,a^*]$ for the edges $a\in Q$, i.e. of
elements $z^{h-2}(c_ae_{head(x)}-c_{a^*} e_{tail(x)})$.
It is clear that these elements span a subspace
of codimension 1 in $A[2(h-2)]$; thus the functional $\Tr$
is unique up to scaling. Moreover, $\Tr(z^{h-2}e_i)$ is
clearly nonzero for any $i$. The proposition is proved.
\end{proof}

Now let $Z$ be the center of $A$.

\begin{corollary}\label{pairi}
The inner product $(x,y)$ defines a nondegenerate pairing
$Z\times A/[A,A]\to \Bbb C$.
\end{corollary}

\begin{proof} The statement is well known but we give a proof for
completeness. If $x\in Z$ and $y=[y_1,y_2]\in [A,A]$ then $(x,y)=
\Tr(x[y_1,y_2])=\Tr([xy_1,y_2])=0$. Thus the pairing in question is
well defined. To show that it is nondegenerate, by Proposition
\ref{trace} (ii), it suffices to
show that $Z^\perp\subset [A,A]$, or equivalently,
$Z\supset [A,A]^\perp$.

The latter statement is obvious. Indeed, if
$\Tr(x[y_1,y_2])=0$ for any $y_1,y_2$, then
$\Tr([x,y_1]y_2)=0$ for any $y_1,y_2$, and therefore
$[x,y_1]=$ for all $y_1$, implying $x\in Z$.
\end{proof}

Let $p(t)=\sum \dim(A/[A,A])[m]t^m$ be the Hilbert polynomial of $A/[A,A]$, and
$p_*(t)=\sum \dim
Z[m]t^m$ be the Hilbert polynomial of $Z$.

\begin{corollary} The polynomials $p,p_*$ are palindromes of each
other, i.e. $p_*(t)=t^{2(h-2)}p(1/t)$.
\end{corollary}

\subsection{The spaces $Z$ and $A/[A,A]$ as $\Bbb C[z]$-modules}

Let $E$ be the subspace of $A$ spanned by
elements $z^je_i$. Obviously, it has dimension $(h-1)r$.
The Hilbert polynomial of $E$ is $\frac{1-t^{2h-2}}{1-t^2}r$.

\begin{proposition}\label{sur} The natural map $\psi: E\to A/[A,A]$ is
surjective.
\end{proposition}

\begin{proof} It is shown in \cite{MOV}, Section 4, that $A_0/[A_0,A_0]$ is
freely spanned by the idempotents $e_i$. This implies that if
$x\in A$ is an element of positive degree $d$
then there exists a homogeneous element $y\in A$ of degree $d-2$
such that $x-zy\in [A,A]$. Thus the statement follows by induction
in $d$.
\end{proof}

Note now that $A/[A,A]$ and $Z$ are naturally $\Bbb C[z]$-modules, and
the pairing $(,)$ between them is invariant in the sense that the
operator of multiplication by $z$ is selfadjoint.

\begin{corollary} \label{gene} The
$\Bbb C[z]$-module $A/[A,A]$ is minimally generated by $e_i$.
\end{corollary}

\begin{proof} Indeed, $A/[A,A]$ is a quotient of $E$ by the
submodule $E\cap [A,A]$, which shows that it is generated by
$e_i$. The minimality of this set of generators is obvious.
\end{proof}

Thus we see that the operator $z$ in $A/[A,A]$ and $Z$
is a direct sum of $r$ nilpotent Jordan blocks, of some sizes
$m_1\le m_2\le...\le m_r$, and
$p(t)=\sum_{i=1}^r\frac{1-t^{2m_i}}{1-t^2}$.

\section{The main theorem}

Let $N_p$ be the number of positive roots for $Q$
of height $p+1$.

One of our main results is the following theorem.

\begin{theorem}\label{flatness} (i) $\dim Z=\dim(A/[A,A])=N$.

(ii) The sizes $m_i$ of the Jordan blocks of $z$ on $Z$ and
$A/[A,A]$ are the exponents of the root system attached to $Q$.
In other words, we have $\dim(A/[A,A])[2p]=N_p$ for all $p\ge
0$.
\end{theorem}

The proof of Theorem \ref{flatness} is given in the next
two subsections.

\subsection{The lower bound}

\begin{proposition}\label{ge}
$\dim Z\ge N$.
\end{proposition}

\begin{proof} According to \cite{ER},
for generic $\lambda$ the algebra $A(\lambda)$ is
semisimple with irreducible representations
$V_\alpha$ corresponding to positive roots $\alpha$.
This implies that the center $Z(\lambda)$
of $A(\lambda)$ is a semisimple algebra of dimension $N$.
Since ${\rm gr}(A(\lambda))=A$, we have ${\rm
gr}(Z(\lambda))\subset Z$, and we get the
desired inequality.
\end{proof}

\subsection{The upper bound}

We have $\sum_p N_p=N$. Theorefore, by Proposition \ref{ge},
to prove Theorem \ref{flatness}, it
suffices to show that $\dim(A/[A,A])[2p]\le N_p$ for all $p\ge
0$.

We do it case by case,
following the idea of
the argument of \cite{MOV}, Section 4.
Since we need to establish the result for generic
$\mu$, it suffices to consider the case $\mu=\rho$.

{\bf Case 1: type $A_n$.} In this case $N_p={\rm max}(n-p,0)$.
Denote the corresponding algebra $A$
by $A^n$, and let us prove the desired statement by induction in
$n$.

The base of induction ($n=1$) is obvious, so let us perform the induction
step. Assume that the statement is known for $A^{n-1}$, and let
us prove it for $A^n$ ($n\ge 2$).

Let $J=Ae_nA$ be the ideal in $A^n$ spanned by
paths passing through the end-vertex $n$ of the
Dynkin diagram. Then $A^n/J=A^{n-1}$. Thus, using the induction
assumption, we see that
$$
\dim A^n/(J+[A^n,A^n])[2p]\le {\rm max}(n-1-p,0).
$$
Thus to establish the induction step (i.e. to show that
$\dim(A^n/[A^n,A^n])\le n-p$),
it suffices to show that $\dim(J/(J\cap [A^n,A^n]))[2p]\le 1$
for $p\le n-1$ and is zero if $p\ge n$.

Define the algebra $B_n:=e_nA^ne_n$.
It is easy to see that the natural map $\phi: B_n\to J/(J\cap
[A^n,A^n])$ is surjective. On the other hand, it is easy to check
that $B_n$ is a commutative $n$-dimensional
algebra: $B_n=\Bbb C[z]/(z^n)$. Thus the desired statement follows.

{\bf Case 2: types $D,E$.} Let $*$ be the nodal vertex of the
Dynkin diagram, and $J=Ae_*A$. Then $A/J=A^{\ell_1}\oplus
A^{\ell_2}\oplus A^{\ell_3}$, where $\ell_1,\ell_2,\ell_3$ are
the lengths of the three legs of the Dynkin diagram.
Thus by Case 1,
$$
\dim A/(J+[A,A])[2p]\le
\sum_{j=1}^3\max(\ell_j-p,0).
$$
So it suffices to show that
$$
\dim J/(J\cap [A,A])[2p]\le N_p',
$$
where $N_p':=N_p-\sum_{j=1}^3\max(\ell_j-p,0)$ is the number
of positive roots of height $p+1$ which contain
the simple root $\alpha_*$ in their expansion.

Define the algebra $B:=e_*Ae_*$. It is easy to see
that the natural map $\phi: B\to J/(J\cap
[A,A])$ is surjective. Thus, it suffices to show that
$$
\dim (B/[B,B])[2p]\le N_p'.
$$

According to \cite{ER}, the algebra $B$ is generated by degree 2
elements $U_1,U_2,U_3$ with defining relations
\begin{equation}\label{eq}
U_1+U_2+U_3=z,\ [z,U_i]=0,
\prod_{m=0}^{\ell_i} (U_i+mz)=0, \ i=1,2,3.
\end{equation}

{\bf Case 2a.} Type $D_{n+2}$ ($n\ge 2$). We have $\ell_1=\ell_2=1$,
$\ell_3=n-1$. So, setting $a=U_1+z/2$, $b=U_2+z/2$, we have the
following defining relations for $B$:
$$
a^2=b^2=z^2/4, [a,z]=[b,z]=0,
$$
$$
(a+b-2z)(a+b-3z)...(a+b-(n+1)z)=0.
$$

Let $a_s=aba...,b_s=bab...$ (words of length $s$).

\begin{lemma}\label{bas}
If $p<n$ then a basis of $B[2p]$ is formed
by the words $a_sz^{p-s}$, $b_sz^{p-s}$, $p\ge s>0$, and $z^p$.
If $p=n$, then a basis of $B[2p]$ is formed
by the words $a_sz^{p-s}$, $b_sz^{p-s}$, $n>s>0$, $a_n$ and
$z^n$. If $p>n$, then a basis of $B[2p]$ is formed
by the words $a_sz^{p-s}$, $b_sz^{p-s}$, $2n-p\ge s>0$, and
$z^p$.
\end{lemma}

\begin{proof}
It is easy to see from the relations that these words are a
spanning set for $B[2p]$. The fact that they are linearly
independent follows from the Hilbert series formula for $B$ given
in \cite{ER}.
\end{proof}


\begin{lemma}\label{lowbound}
One has $\dim (B/[B,B])[2p]\le N_p'$.
\end{lemma}

\begin{proof}
It is straightforward to show (by explicit inspection of the root
system of type $D_{n+2}$) that $N_0'=1$, $N_p'=3+[p/2]$ if $1\le
p\le n-1$, $N_n'=2+[n/2]$, and $N_p=1+[n-p/2]$ for $p>n$, where
$[x]$ is the integer part of $x$. For odd $s>1$ and $p\ge s$, we
have
$$a_sz^{p-s}=\frac{1}{4}b_{s-2}z^{p-s+2}+[a,b_{s-1}z^{p-s}],$$
$$b_sz^{p-s}=\frac{1}{4}a_{s-2}z^{p-s+2}+[b,a_{s-1}z^{p-s}].$$
Also, for even $s>0$,
$$(a_s-b_s)z^{p-s}=[a,b_{s-1}z^{p-s}].$$
This together with Lemma \ref{bas} implies that $(B/[B,B])[2p]$ is
spanned by $z^p$, $az^{p-1}$, $bz^{p-1}$, and $a_sz^{p-s}$ for
even $s>0$. Hence, for $p\ge 1$ we have
$$\dim (B/[B,B])[2p]\le 3+[p/2],$$
i.e. the Lemma is proved for $p<n$. Moreover, for $p=n$ the last
relation of $B$ implies that $z^n$ is a linear combination of
$a_sz^{n-s}$ and $b_sz^{n-s}$ for $s>0$, which implies that
$$\dim (B/[B,B])[2n]\le 2+[n/2],$$ i.e. the Lemma is also proved for $p=n$.

Now let us prove the lemma for $p>n$. Let $Z_B$ be the center of
$B$. The pairing $(x,y)$ of Proposition \ref{trace} has degree
$4n$. Therefore, by Proposition \ref{trace} (similarly to
Corollary \ref{pairi}), it suffices to show that $\dim Z_B[2k]\le
1+[k/2]$ for $k<n$. In showing this, we can obviously ignore the
last relation of $B$ (which has degree $2n$). In other words, we
should consider the algebra $B'$ with generators $a,b,z$ and
relations $a^2=b^2=z^2/4$, $[z,a]=[z,b]=0$. It is easy to see that
a basis in $B'$ is formed by elements $z^p(a+b)^{2q}$,
$az^p(a+b)^{2q}$, $bz^p(a+b)^{2q}$, $abz^p(a+b)^{2q}$, $p,q\ge 0$,
and thus the center of $B'$ is spanned by $z^p(a+b)^{2q}$, which
implies the desired inequality. \end{proof}

{\bf Case 2b.} Types $E_6,E_7,E_8$.
Using the presentation (\ref{eq}) of $B$ and the Magma code by the third author
\cite{Mag}, one determines, by a direct computer calculation,
that $\dim (B/[B,B])[2p]=N_p'$.

Theorem \ref{flatness} is proved.

\subsection{Derivations of $A$.}

Theorem \ref{flatness} implies the following result.

\begin{corollary}\label{inn}
Every derivation of $A$ which annihilates $R$ and $z$ is inner.
\end{corollary}

{\bf Remark.} In this corollary, we can omit the hypothesis that the derivation
annihilates $R$. Indeed, for any derivation $D$ of $A$, if we let
$u_D:=\sum_{i\in I} e_i D(e_i)$,
then $D+\ad(u_D)$ annihilates $R$ and has the same action as $D$ on the
central element $z$.

\begin{proof}
We consider the complex of graded vector spaces
$$
0\to D_0\to D_1\to D_2\to 0,
$$
 with
differentials $d_i: D_i\to D_{i+1}$, where
$D_0=A^R[2]$, $D_1=(A\otimes_R V)^R$, $D_2=A^R$ (where the superscript $R$
denotes the $R$-invariants in a bimodule, and $[i]$ denotes the
shift of grading), and
$$
d_0(x)=\sum_{a\in Q}([x,a]\otimes a^*-[x,a^*]\otimes a),
$$
$$
d_1(y\otimes b)=[y,b].
$$
It is clear that these differentials have degree $0$.
The fact that $d_1\circ d_0=0$
follows from the fact that $z$ is a central element.

Let $H_0,H_1,H_2$ be the homology groups of the complex
$D_\bullet$. Then we have $H_0=Z[2]$, $H_2=A/[A,A]$.

Let $q(t)$ be the Hilbert polynomial of $H_1$.
Then, computing the Euler characteristic
in each homogeneous component of $D_\bullet$,
we obtain the following identity for Hilbert polynomials:
$$
t^2p_*(t)+p(t)-q(t)=\Tr((1-Ct+t^2)h(t)),
$$
where $h(t)$ is the (matrix valued) Hilbert polynomial of $A$,
and $C$ is the adjacency matrix of $\overline{Q}$.
But it is proved in \cite{ER} that
$$
h(t)=\frac{1-t^{2h}}{1-t^2}(1-Ct+t^2)^{-1}.
$$
This implies that
$$
q(t)=t^2p_*(t)+p(t)-\frac{1-t^{2h}}{1-t^2}r.
$$
Now recall that the exponents of a root system satisfy the
equality $m_{r+1-i}=h-m_i$. This implies that
$t^2p_*(t)+p(t)=\frac{1-t^{2h}}{1-t^2}r$, and hence $q(t)=0$.
Thus $H_1=0$.

Now let $D$ be a derivation of $A$ which annihilates $R$ and $z$.
Let $x_D:=\sum_{a\in Q}(Da\otimes a^*-Da^*\otimes a)$.
Then $d_1x_D=0$. Since $H_1=0$, this implies that $x_D=d_0y$,
i.e. $D={\rm ad}y$, as desired. The corollary is proved.
\end{proof}

\section{Relation to simple Lie algebras}

The computer assisted case-by-case proof of Theorem \ref{flatness} makes it
look mysterious (especially part (ii)). The results
of this section demystify this theorem, by making explicit the
relation of the structure of $A/[A,A]$ with that of the
maximal nilpotent subalgebra of the simple Lie algebra
corresponding to $Q$.

\subsection{The results}

Let us color the vertices of $Q$ white and black so that
every edge connects a white vertex with a black vertex.
Let $\varepsilon_i$ be $+1$ for white vertices $i$ and $-1$ for
black vertices. Let $F=\sum_i \varepsilon_i F_i$ be a principal
nilpotent element.

Let $h_\lambda\in {\mathfrak{h}}$ be the
element corresponding to the weight $\lambda\in \h^*$ under the standard inner
product on ${\mathfrak{h}}^*$ normalized so that $(\alpha,\alpha)=2$ for
roots $\alpha$.

The following theorem characterizes explicitly the space $E\cap
[A,A]$.

\begin{theorem}\label{ff}
Let $\phi_i, i\in I$ be complex numbers, and $s\ge 0$
be an integer. Then the element
$z^s(\sum_i \varepsilon_i\phi_i e_i)$ is in $[A,A]$ if and only if
$$
({\rm ad}(F){\rm ad}(h_\mu)^{-1})^s(\sum \phi_i F_i)=0
$$
in $\n$.
\end{theorem}

Note that Theorem \ref{ff} implies Theorem \ref{flatness}.
In the proof of Theorem \ref{ff}, given in the next subsection,
we will use only part (i)
of Theorem \ref{flatness}, so we obtain a new proof of Theorem
\ref{flatness}, part (ii).

The result of Theorem \ref{ff} can be stated more explicitly as
follows.

Let $V_i$ be the space of complex-valued functions on the
set of positive roots for the quiver $Q$ of height $i$
(i.e. sums of $i$ simple roots).
Define the operator
$T_i: V_i\to V_{i+1}$ by the formula
$$
(T_if)(\alpha)=\sum_{j:(\alpha_j,\alpha)=1}
\frac{f(\alpha-\alpha_j)}{(\mu,\alpha-\alpha_i)}.
$$
(Note that $\alpha-\alpha_j$ is a root iff
$(\alpha,\alpha_j)=1$.)

\begin{theorem}\label{pt} Let $\phi\in V_1$,
$\phi_i=\phi(\alpha_i)$. Let $s\ge 0$. Then element
$z^s(\sum_i \varepsilon_i \phi_ie_i)$ is in $[A,A]$ iff
$T_sT_{s-1}...T_1\phi=0$.
\end{theorem}

\begin{proof}
According to \cite{Lu}, there is a Chevalley basis
$\lbrace{F_\alpha\rbrace}$ of $\n$ normalized in such a way that
$[F_i,F_\alpha]=\varepsilon_i F_{\alpha+\alpha_i}$ provided
$\alpha+\alpha_i$ is a root.
Therefore,
$$
{\rm ad}(F){\rm ad}(h_\mu)^{-1}
\sum_{\beta\in {\mathcal R}: (\rho,\beta)=d}\phi_\beta
F_\beta=
$$
$$
\sum_{\gamma\in {\mathcal R}:
(\rho,\gamma)=d+1}
\sum_i\frac{\phi_{\gamma-\alpha_i}}{(\mu,\gamma-\alpha_i)} F_\gamma.
$$
Thus Theorem \ref{ff} implies
Theorem \ref{pt}.
\end{proof}

\begin{corollary} The explicit form of the trace functional for
$A$ is
$$
\Tr(z^{h-2}e_i)=\varepsilon_i T_{h-2}...T_1u_i,
$$
where $u_i\in V_1$ is such that $u_i(\alpha_j)=\delta_{ij}$.
\end{corollary}

This formula can be written more explicitly as follows.
Let $\theta$ be the maximal root of ${\mathcal R}$.
Define a path in ${\mathcal R}$ to be a sequence of positive roots
$\beta_1,\beta_2,...,\beta_m$ such that
$\beta_{i+1}-\beta_i=\alpha_{j_i}$ for some $j_i$.
Define weight of such a path $\pi$ to be
$$
w_\mu(\pi)=\prod_{i=1}^{m-1}(\mu,\beta_i)^{-1}.
$$
Then we get
$$
\Tr(z^{h-2}e_i)=\varepsilon_i\sum_\pi w_\mu(\pi),
$$
where the summation is taken over all paths $\pi$
which start at $\alpha_i$ and end at $\theta$ (so they have
length $h-1$). In particular, if $\mu=\rho$ then after
renormalization we get
$$
\Tr(z^{h-2}e_i)=\varepsilon_in_i,
$$
where $n_i$ is the number of paths leading from $\alpha_i$ to
$\theta$.

\subsection{Proof of Theorem \ref{ff}}

Let $W(\lambda)$ be the space of collections of polynomials
$f_i,i\in I$ of degree $\le h-2$,
such that $\sum_{i\in I} f_i(z)\varepsilon_i e_i\in [A(\lambda),A(\lambda)]$.

\begin{proposition}\label{lamb} Let $\lambda$ be generic.
Let $f_i$, $i\in I$, be polynomials of degree $\le h-2$.
Then $\lbrace{f_i,i\in I\rbrace}$ belongs to $W(\lambda)$ iff
$$
\sum
f_i\biggl(-\frac{(\lambda,\alpha)}{(\mu,\alpha)}\biggr)
\varepsilon_i (\alpha,\omega_i)=0
$$
for all positive roots $\alpha$.
\end{proposition}

\begin{proof} Let us calculate the trace
of $\sum f_i(z)\varepsilon_i e_i$ in the irreducible
representation $V_\alpha$ of $A(\lambda)$ whose
dimension vector is $\alpha$. Since $z$ acts on this
representation by the scalar
$-\frac{(\lambda,\alpha)}{(\mu,\alpha)}$,
we get the statement.
\end{proof}

\begin{proposition}\label{lamb1}
Let $\lambda$ be generic.
Let $f_i$, $i\in I$, be polynomials of degree $\le h-2$.
Then $\lbrace{f_i,i\in I\rbrace}$ belongs to $W(\lambda)$ iff
$$
\sum_{i\in I} f_i({\rm ad}(-h_\lambda+F){\rm ad}(h_\mu)^{-1})F_i=0.
$$
\end{proposition}

\begin{proof}
The linear operator $L:={\rm ad}(-h_\lambda+F){\rm ad}(h_\mu)^{-1}$ on $\n$
has eigenvalues
$-(\lambda,\alpha)/(\mu,\alpha)$, where $\alpha$ ranges over positive roots;
let the corresponding eigenvectors be $v_\alpha$.  Then, if we write
\[
F_i = \sum_\alpha c_i(\alpha) v_\alpha,
\]
we find
\[
\sum_{i\in I} f_i({\rm ad}(-h_\lambda+F){\rm ad}(h_\mu)^{-1})F_i
=
\sum_{i\in I} \sum_\alpha f_i\lbrace(-\frac{(\lambda,\alpha)}
{(\mu,\alpha)}\rbrace) c_i(\alpha)
v_\alpha.
\]
In particular, to prove the proposition, it will suffice to show that
\[
c_i(\alpha) \propto \varepsilon_i(\alpha,\omega_i).
\]

Now, by duality, $c_i(\alpha)$
can be computed as the coefficient of
$E_i\in \n_+$ in the expansion of the eigenvector $v_\alpha^*$
of the dual operator
\[
L^*={\rm ad}(h_\mu)^{-1}{\rm ad}(-h_\lambda+F)
\]
on ${\mathfrak g}$ with eigenvalue $-(\lambda,\alpha)/(\mu,\alpha)$.

Let $X_\alpha\in \n_+$ be the projection of $v_\alpha^*$
to $\n_+$ along $\h\oplus\n_-$. Then the element
\[
y:=[-h_\lambda+F+\frac{(\lambda,\alpha)}{(\mu,\alpha)} h_\mu,X_\alpha]
\]
must belong to the Cartan subalgebra $\h$.

Recall that $\lambda\in \h^*$ is generic.
Therefore, if $\nu\in \h^*$ is any element of the orthogonal
complement of $\alpha$, then there exists
${\mathcal N}\in {\frak n}$ such that
\[
[-h_\lambda+F+\frac{(\lambda,\alpha)}{(\mu,\alpha)} h_\mu,h_\nu + {\mathcal N}] = 0,
\]
and thus $(y,h_\nu+{\mathcal N}) = 0$.  It follows that $y\propto h_\alpha$; since
$X_\alpha$ was only determined up to scale, we may as well insist that
$y=h_\alpha$.

Since $X_\alpha = \sum c_i(\alpha) E_i + \text{lower terms}$,
we find that
\[
h_\alpha = [F,X_\alpha]_\h = \sum \varepsilon_i c_i(\alpha) h_{\alpha_i}
\]
(where the subscript $\h$ denotes the $\h$-part),
and thus $\varepsilon_i c_i(\alpha) = (\alpha,\omega_i)$ as required.
\end{proof}

Now we can finish the proof of Theorem \ref{ff}.
For this, note that by Theorem \ref{flatness}(i),
the space $W(0)$ is the limit of spaces $W(\lambda)$ as
$\lambda\to 0$. Therefore, Proposition \ref{lamb1}
implies Theorem \ref{ff}.

\end{document}